\documentclass[aip,reprint,groupedaddress]{revtex4-1}
\usepackage{color}
\usepackage{graphicx}
\usepackage{bm}
\draft
\begin{document}
\title{Simulated bifurcation for higher-order cost functions}
\author{Taro Kanao}
\email[]{taro.kanao@toshiba.co.jp}
\author{Hayato Goto}
\affiliation{Frontier Research Laboratory, Corporate Research \& Development Center, Toshiba Corporation, 1, Komukai-Toshiba-cho, Saiwai-ku, Kawasaki 212-8582, Japan}
\date{\today}
\begin{abstract}
High-performance Ising machines for solving combinatorial optimization problems have been developed with digital processors implementing heuristic algorithms such as simulated bifurcation (SB). Although Ising machines have been designed for second-order cost functions, there are practical problems expressed naturally by higher-order cost functions. In this work, we extend SB to such higher-order cost functions. By solving a problem having third-order cost functions, we show that the higher-order SB can outperform not only the second-order SB with additional spin variables, but also simulated annealing applied directly to the third-order cost functions. This result suggests that the higher-order SB can be practically useful.
\end{abstract}
\maketitle

Special purpose hardware solvers for the Ising problem~\cite{Barahona1982} have been applied to combinatorial optimization problems~\cite{Lucas2014}.
Such hardware solvers are referred to as Ising machines~\cite{Mohseni2022a}.
Ising machines with analogue computers have attracted attention, such as quantum annealers~\cite{Kadowaki1998} with superconducting circuits~\cite{Johnson2011} and coherent Ising machines with pulse lasers~\cite{Marandi2014}.
Also, Ising machines with digital processors have been developed based on, e.g., simulated annealing (SA)~\cite{Kirkpatrick1983, Yamaoka2016, Aramon2019}, simulated bifurcation (SB)~\cite{Goto2019}, and other algorithms~\cite{Leleu2021a}.

SB is a heuristic algorithm based on numerical simulation of a classical counterpart of a quantum bifurcation machine, which is a quantum annealer using nonlinear oscillators~\cite{Goto2016, Puri2017, Goto2019a, Kanao2021}.
SB allows us to simultaneously update variables, and thus can easily be accelerated by parallel processing with, e.g., graphics processing units (GPUs) and field-programmable gate arrays (FPGAs)~\cite{Goto2019, Tatsumura2019, Zou2020, Tatsumura2021}.
Two improved variants of SB have been reported: ballistic and discrete SBs (bSB and dSB)~\cite{Goto2021}.
In addition, bSB and dSB with thermal fluctuations have been studied~\cite{Kanao2022a}.
SB has been benchmarked against other Ising machines~\cite{Oshiyama2022, Kowalsky2022}.
Besides, SB has been applied to financial portfolio optimizations~\cite{Steinhauer2020a, Cohen2020}, financial trading machines~\cite{Tatsumura2020}, traveling salesman problems~\cite{Zhang2021c}, clustering with a hybrid method~\cite{Matsumoto2022}, and wireless communications~\cite{Zhang2022c}.

In the Ising problem, the aim is to minimize cost functions written in second-order polynomials of spin variables, for which Ising machines have been designed.
However, some combinatorial optimization problems are naturally expressed by cost functions with higher-order polynomials than second order (higher-order cost functions).
To apply Ising machines to such problems, higher-order cost functions have been transformed to second-order ones~\cite{Kolmogorov2004, Ishikawa2011, Xia2018, Dattani2019, Fujisaki2022a}, though such transformations usually require additional spin variables and thus increase computational costs.

In this work, we extend SB to apply it directly to a problem with higher-order cost functions.
Since SB is based on simulation, this approach can easily be implemented.
To benchmark the higher-order SB, we use a three-regular three-XORSAT (3R3X) problem~\cite{Hen2019, Bernaschi2021, Bellitti2021, Kowalsky2022}, because this problem yields third-order cost functions and its solution can be known in advance.
We compare the third-order SB with that transformed to second order.
As a result, we find that third-order bSB and dSB (3bSB and 3dSB) show higher performance than the second-order ones (2bSB and 2dSB).
Surprisingly, 3bSB performs better than 3dSB, while 2dSB does better than 2bSB as in the previous study~\cite{Goto2021}.
We further find that 3bSB can be superior to SA applied directly to the third-order cost functions.
These results suggest that the higher-order SB will offer  another useful approach to special purpose hardware solvers for combinatorial optimization problems.

A third-order cost function with $N$ spin variables $s_i=\pm1$ can generally be written as
\begin{eqnarray}
	E(\bm{s})&=&-\!\sum_i\!h_is_i-\!\sum_{i,j}\!J_{ij}s_is_j-\!\sum_{i,j,k}\!K_{ijk}s_is_js_k,\label{eq_cost}
\end{eqnarray}
where $\bm{s}$ is a vector consisting of $s_i$, and $h_i, J_{ij}$, and $K_{ijk}$ are coefficients.
Since $s_i^2=1$, we assume without loss of generality that $J_{ii}=0$, and $K_{ijk}=0$ if two of $i, j$, and $k$ are equal.

We here extend second-order bSB and dSB to third-order ones.
The third-order bSB and dSB are based on the following equations of motion for positions $x_i$ and momenta $y_i$ corresponding to $s_i$:
\begin{eqnarray}
	\dot{y}_i&=&-\left[a_0-a(t)\right]x_i+cf_i,\label{eq_y}\\
	\dot{x}_i&=&a_0y_i,\label{eq_x}
\end{eqnarray}
where the dots denote time derivatives, $a_0$ is a constant parameter, $a(t)$ is a bifurcation parameter increased from 0 to $a_0$ with the time $t$, and $c$ is a normalization factor determined from forces $f_i$.
The forces $f_i$ are given by
\begin{eqnarray}
	f_i&=&-\frac{\partial E(\bm{x})}{\partial x_i}=G_i(\bm{x})\quad\text{(bSB)},\label{eq_f_b}\\
	f_i&=&-\left.\frac{\partial E(\bm{x})}{\partial x_i}\right|_{\bm{x}=\bm{s}}=G_i(\bm{s})\quad\text{(dSB)}.\label{eq_f_d}
\end{eqnarray}
For 3bSB, $f_i$ are the gradients $G_i(\bm{x})$ of $E(\bm{x})$, while for 3dSB, $\bm{x}$ in the gradients is substituted by $\bm{s}$, which is a vector of the signs of $x_i$, $s_i={\rm sgn}\left(x_i\right)$.
Both 3bSB and 3dSB assume perfectly inelastic walls at $x_i=\pm1$~\cite{Goto2021}, that is, if $\left|x_i\right|>1$, $x_i$ and $y_i$ are set to ${\rm sgn}\left(x_i\right)$ and $0$, respectively.
Equations~(\ref{eq_y}) and (\ref{eq_x}) are solved with the symplectic Euler method~\cite{Leimkuhler2004}, which discretizes the time with an interval $\Delta t$ and computes time evolutions of $x_i$ and $y_i$ step by step.
The signs of $x_i$ at the final time give a solution.
The solution at least corresponds to a local minimum of $E(\bm{s})$ as for 2bSB and 2dSB~\cite{Goto2021}, because Eq.~(\ref{eq_y}) becomes $\dot{y}_i=c_0f_i$ at the final time.
Note that when all $K_{ijk}$ are zero in Eq.~(\ref{eq_cost}), Eqs.~(\ref{eq_y})-(\ref{eq_f_d}) reproduce 2bSB and 2dSB in Ref.~\onlinecite{Goto2021}.
We thus focus on the third term in Eq.~(\ref{eq_cost}) in the following.

To reduce the computational cost of calculating $f_i$, we utilize sparse connectivity of the 3R3X problem, that is, the number, $M$, of nonzero $K_{ijk}$ is much smaller than the total number of $K_{ijk}$.
We expect that sparse connectivities often appear in applications and therefore the following method will be useful.
By denoting nonzero $K_{ijk}$ by $K_m$ $(m=1, 2, \cdots, M)$, $E(\bm{s})$ is expressed as
\begin{eqnarray}
	E(\bm{s})&=&-\sum_{m=1}^MK_ms_{v_{m1}}s_{v_{m2}}s_{v_{m3}},\label{eq_sparse}
\end{eqnarray}
where $v_{mn}$ with $n=1, 2, 3$ represent the indices of spin variables included in the $m$th term.
To further reduce the time required for evaluating the third-order SB, we calculate the gradients by
\begin{eqnarray}
	G_i(\bm{x})&=&\frac{1}{x_i+\epsilon}\sum_{m=1}^Ma_{im}K_mx_{v_{m1}}x_{v_{m2}}x_{v_{m3}}\quad\text{(bSB)},\label{eq_def_b}\\
	G_i(\bm{s})&=&s_i\sum_{m=1}^Ma_{im}K_ms_{v_{m1}}s_{v_{m2}}s_{v_{m3}}\quad\text{(dSB)},\label{eq_def_d}\\
	a_{im}&=&\delta_{i,v_{m1}}+\delta_{i,v_{m2}}+\delta_{i,v_{m3}},
\end{eqnarray}
where $\delta_{i,j}$ is 1 if $i=j$ and otherwise 0.
Equation~(\ref{eq_def_d}) exactly gives the derivative in Eq.~(\ref{eq_f_d}), because $a_{im}$ pick up terms including $s_i$ from Eq.~(\ref{eq_sparse}), and these $s_i$ are canceled by multiplying $s_i$ in Eq.~(\ref{eq_def_d}) with $s_i^2=1$.
This method is valid since Eq.~(\ref{eq_sparse}) is linear with respect to each $s_i$.
Similarly, Eq.~(\ref{eq_def_b}) yields the same value as Eq.~(\ref{eq_f_b}) if $x_i\neq0$ and $\epsilon=0$.
To stabilize calculations, we set $\epsilon=10^{-14}$.
We empirically found that Eqs.~(\ref{eq_def_b}) and (\ref{eq_def_d}) could reduce the evaluation time compared with direct differentiation, which may be because these equations gather three terms derived from $x_{v_{m1}}x_{v_{m2}}x_{v_{m3}}$ together and thus reduce memory accesses for $x_{v_{mn}}$ and $s_{v_{mn}}$.

The 3R3X problem is mapped to a minimization problem of Eq.~(\ref{eq_sparse}) with the following $v_{mn}$ and $K_m$~\cite{Hen2019}.
Here, $M$ equals $N$.
For each of $n=1, 2, 3$, the vector of indices $v_{mn}$ is a permutation of $1, 2, \cdots, N$.
For each $m$, three indices $v_{m1}$, $v_{m2}$, and $v_{m3}$ are different from each other.
For different $m$, the combinations of $v_{m1}$, $v_{m2}$, and $v_{m3}$ are different.
$K_m$ is given by $(-1)^{b_m}$ with $b_m=0, 1$.
Thus $E(\bm{s})\geq-N$ holds.

For a 3R3X instance, a solution $\bm{s}_0$ minimizing $E(\bm{s})$ can be obtained with linear algebra~\cite{Hen2019}, that is, $\bm{s}_0$ is related to a solution, $\xi_i=0, 1$, of linear equations,
\begin{eqnarray}
	\sum_{i=1}^Na_{mi}\xi_i&=&b_m,\quad\mathrm{mod}\quad 2,\label{eq_linear}
\end{eqnarray}
by $s_{i0}=(-1)^{\xi_i}$.
However, depending on $a_{mi}$ and $b_m$, Eq.~(\ref{eq_linear}) may have no solution.
To ensure that Eq.~(\ref{eq_linear}) has at least one solution, we instead give $a_{mi}$ and $\xi_i$, and determine $b_m$ with Eq.~(\ref{eq_linear}).
Then an exact solution $\bm{s}_0$ is indicated by $E(\bm{s}_0)=-N$.

A third-order cost function of the 3R3X problem can be decomposed to second-order one by mapping a term  $-(-1)^bs_1s_2s_3$ to
\begin{eqnarray}
&&h\left(s_1+s_2+s_3\right)+\tilde{h}\tilde{s}\nonumber\\
	&&+J\left(s_1s_2+s_2s_3+s_3s_1\right)+\tilde{J}\left(s_1+s_2+s_3\right)\tilde{s},\label{eq_gadget}
\end{eqnarray}
where $\tilde{s}=\pm1$ is an ancillary spin variable, and $J$, $\tilde{J}$, $h$, and $\tilde{h}$ are coefficients determined so that $s_1$, $s_2$, and $s_3$ minimizing Eq.~(\ref{eq_gadget}) coincide with ones minimizing $-(-1)^bs_1s_2s_3$~\cite{Kolmogorov2004, Leib2016, Hen2019}.
We use $J=1/4$, $\tilde{J}=(-1)^b/2$, $h=-(-1)^b/4$, and $\tilde{h}=-1/2$, which normalize the minimum of Eq.~(\ref{eq_gadget}) to $-1$.
The ancillary spin variables increase the total number of spin variables to $2N$.

We also directly apply SA to third-order cost functions (3SA).
In an update in SA~\cite{Kirkpatrick1983}, one calculates a change in a cost function $\Delta E_i$ caused by a flip of $s_i$.
The flip is accepted if $\beta\Delta E_i<-\ln R$, where $\beta$ and $R$ are the inverse temperature increased during annealing and a random number in the interval $(0,1)$, respectively.
In a step of SA, we sequentially update spin variables $s_i$~\cite{Isakov2015}.

In numerical experiments, we measure computational costs of finding exact solutions as follows.
We randomly generate 100 3R3X instances for each $N$.
For these instances, we run the above mentioned SB and SA many times with different initial conditions or random numbers, and obtain probabilities of finding exact solutions, $P$.
Then, using the number of steps in a run, $N_{\mathrm{s}}$, we estimate the number of steps necessary to obtain exact solutions with a probability of 99\% (step-to-solution) by~\cite{Leleu2021a, Kanao2022a}
\begin{eqnarray}
	S&=&N_{\mathrm{s}}\frac{\log0.01}{\log(1-P)}.
\end{eqnarray}
A smaller $S$ indicates that an algorithm provides an exact solution more efficiently (or faster if computation time for each step is assumed to be the same).

We set the parameters as follows.
For SB, $a_0=1$, and $a(t)$ is increased linearly.
To normalize the term $cf_i$ in Eq.~(\ref{eq_y}), we give $c$ by~\cite{Sakai2020}
\begin{eqnarray}
	c&=&c_1\left(\frac{1}{\nu}\sum_{i=1}^{\nu}f_i^2\right)^{-1/2},\label{eq_c0}
\end{eqnarray}
where $c_1$ and $\nu$ are a constant parameter and the number of spin variables, respectively ($\nu$ are $2N$ and $N$ for the second- and third-order SBs, respectively).
We optimize $N_{\mathrm{s}}$ for each $N$, and $(\Delta t, c_1)$ for certain large values of $N$, to minimize the median of $S$ for the 100 instances~\cite{Kowalsky2022}.
Initial $x_i$ and $y_i$ are generated by uniform random numbers in the interval $(-1,1)$.
To suppress a statistical error in $P$, the number of runs $N_{\mathrm{r}}$ is set to be large such that $N_{\mathrm{r}}P\gg1$~\cite{Goto2021}, though this condition is not necessarily met for the largest $N$ in this work owing to limitations of computational resources.
For SA, $\beta$ is linearly increased from zero to a final value $\beta_1$, and $\beta_1$ and $N_{\mathrm{s}}$ are optimized as in SB.

\begin{figure}
	\centering
	\includegraphics[width=8cm]{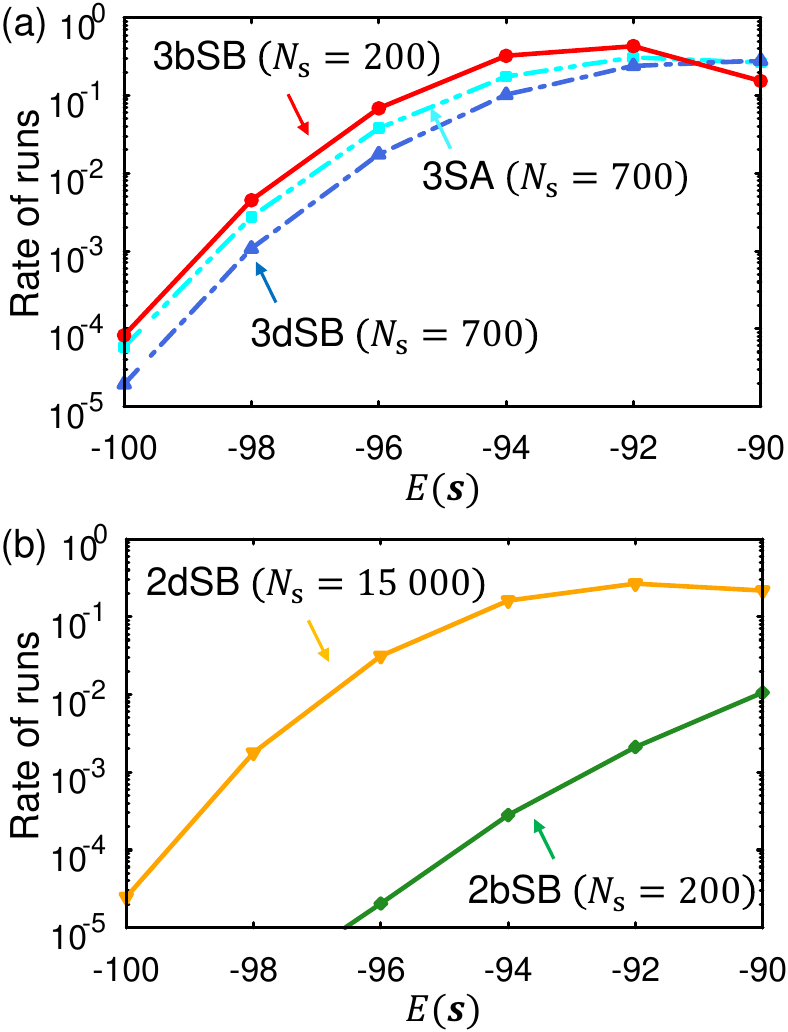}
	\caption{Distributions of a third-order cost function $E(\bm{s})$ for a 3R3X instance with 100 spin variables.
		(a) 3bSB, 3dSB, and 3SA.
		(b) 2bSB and 2dSB.
		The lines are eye guides.
		The numbers of steps $N_{\mathrm{s}}$ are shown in the parentheses.
		Here, $N_{\mathrm{s}}$ and the following parameters are determined to minimize step-to-solution $S$.
		The parameters are $\left(\Delta t, c_1\right)=(1.1, 0.7)$ for 3bSB, $(0.7, 1.1)$ for 3dSB, $(0.8, 0.9)$ for 2bSB, $(0.7, 1.6)$ for 2dSB, and $\beta_1=2$ for 3SA.
		\label{fig_rate_e}}
\end{figure}

We first present typical distributions of $E(\bm{s})$ for a 3R3X instance, for which 3bSB gives $S$ near the median of the 100 instances.
Figure~\ref{fig_rate_e} shows the distributions by the rates of runs resulting in $E(\bm{s})$.
A higher value of the rate means a higher probability of finding solutions with certain accuracy, though here the parameters are optimized to minimize $S$ as mentioned above and the rates can be improved by increasing $N_{\mathrm{s}}$.
Interestingly, Fig.~\ref{fig_rate_e}(a) shows that the rate for small $E(\bm{s})$ is higher for 3bSB than for 3dSB despite smaller $N_{\mathrm{s}}$ for 3bSB.
On the other hand, Fig.~\ref{fig_rate_e}(b) shows that the rate is higher for 2dSB than for 2bSB in accordance with much larger $N_{\mathrm{s}}$ for 2dSB.
As shown in Fig.~\ref{fig_rate_e}(a), the rate for 3SA is between those of 3bSB and 3dSB.
These rates and $N_{\mathrm{s}}$ suggest the highest performance of 3bSB.

\begin{figure}
	\centering
	\includegraphics[width=8.5cm]{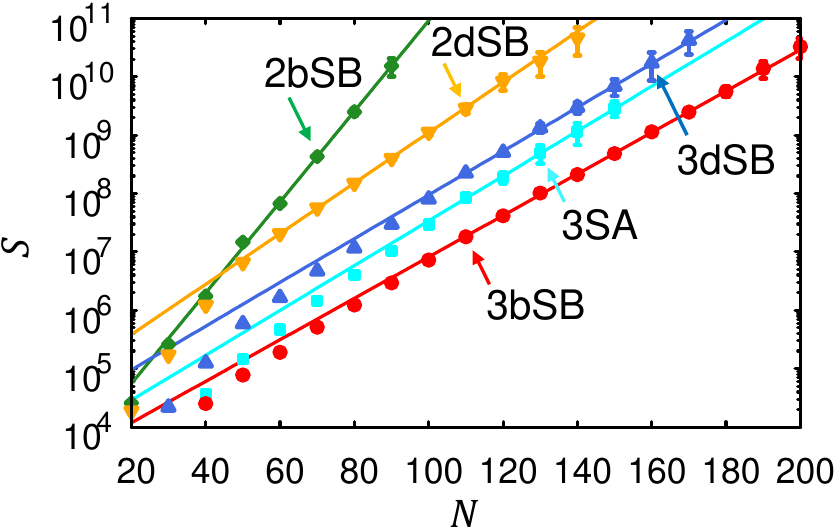}
	\caption{The medians of step-to-solution $S$ for 100 3R3X instances as functions of the number of spin variables $N$.
		The lines and error bars show fittings by Eq.~(\ref{eq_scale}) and statistical errors~\cite{Goto2021}, respectively.
		Parameters $\Delta t, c_1$, and $\beta_1$ are the same as in Fig.~\ref{fig_rate_e}.
		\label{fig_s_n}}
\end{figure}
We next evaluate the results for 100 3R3X instances.
Figure~\ref{fig_s_n} shows the medians of $S$ as functions of $N$.
All the algorithms exhibit exponential increases of $S$ with $N$.
The values of $S$ are smaller for third-order bSB and dSB than for the second-order ones.
Among the second-order ones, 2dSB gives smaller $S$ than 2bSB, which is consistent with a reported result that dSB more often reaches exact solutions for the Ising problem~\cite{Goto2021}.
Among the third-order SB, in contrast, 3bSB shows smaller $S$ than 3dSB.
The higher performance of 3bSB might be because bSB can rapidly converge to a local minimum of a cost function~\cite{Goto2021} and a local search algorithm is effective for the 3R3X problem~\cite{Bernaschi2021, Bellitti2021, Kowalsky2022}.
A comparison between 2bSB and 3bSB suggests that the local search by bSB is especially effective when directly applied to third-order cost functions with fewer spin variables.
3SA results in $S$ between those of 3bSB and 3dSB, indicating that 3bSB is superior to 3SA.

Scaling can be measured by exponents $\alpha$ and $\beta$ such that
\begin{eqnarray}
	S&\simeq&10^{\alpha N+\beta},\label{eq_scale}
\end{eqnarray}
for large $N$~\cite{Kowalsky2022}.
A smaller $\alpha$ allows an algorithm to be applied to larger instances.
We fit the medians of $S$ by Eq.~(\ref{eq_scale}) with the least squares method as the lines in Fig.~\ref{fig_s_n}.
Here, only the data for intermediate $N$ are used for the fittings, because errors due to finite-size effects and statistical errors are large for small and large $N$, respectively~\cite{Kowalsky2022}.
Table~\ref{tab_exponent} shows $\alpha$ and $\beta$.
The value of $\alpha$ for 3bSB is the smallest.
$\alpha$ for 3dSB is a little larger than it, and $\alpha$ for 2dSB is significantly larger.
$\alpha$ for 2bSB is by far the largest.
The value of $\alpha$ for 3SA is larger than those for 3bSB and 3dSB.
These $\alpha$ demonstrate that 3bSB scales the best among these algorithms.

\begin{table}
	\caption{Scaling exponents.
		The numbers in the parentheses are standard deviations in the least squares method.}
	\label{tab_exponent}
	\centering
	\begin{ruledtabular}
		\begin{tabular}{cccccc}
			&3bSB&3dSB&3SA&2bSB&2dSB\\
			\hline
			$\alpha$&0.0355(2)&0.0375(6)&0.0384(6)&0.078(1)&0.0433(2)\\
			$\beta$&3.36(4)&4.23(8)&3.69(8)&3.19(9)&4.72(2)\\
		\end{tabular}
	\end{ruledtabular}
\end{table}

Finally, we compare $\alpha$ with the results of Ref.~\onlinecite{Kowalsky2022}, where various Ising machines are benchmarked using the 3R3X problem expressed by second-order cost functions.
Among the Ising machines, Digital Annealer reported the smallest exponent corresponding to $\alpha=0.0370(8)$~\cite{Footnote2022}.
Compared with this value, $\alpha$ for 3bSB seems to be smaller.
$\alpha$ for Simulated Bifurcation Machine (provided by Toshiba Digital Solutions Corporation~\cite{TDSL}) corresponds to $0.0434(12)$, which coincides with that for 2dSB in this work within the range of errors.
This coincidence may support the validity of our evaluation.
In Ref.~\onlinecite{Kowalsky2022}, a quasi-greedy algorithm called SATonGPU~\cite{Bernaschi2021, Bellitti2021, Kowalsky2022} is also compared, which is specialized for the 3R3X problem and deals with third-order cost functions directly.
As a result, SATonGPU provided the exponent corresponding to $\alpha=0.0341(2)$~\cite{Footnote2021}, which is close to an estimated optimal exponent~\cite{Bernaschi2021}.
The value of $\alpha$ for SATonGPU is smaller than $\alpha$ for 3bSB in our work, but the difference is fairly small.
This implies that 3bSB can achieve almost optimal scaling even though 3bSB is not specialized for the 3R3X problem.

In summary, we have demonstrated that third-order SB can perform better than second-order one for a 3R3X problem.
Although 2dSB provides better performance than 2bSB, 3bSB instead does better than 3dSB.
Furthermore, we have shown that 3bSB can be superior to 3SA.
These results indicate that higher-order SB can be useful for solving combinatorial optimization problems expressed naturally by higher-order cost functions.
Since we have formulated the higher-order SB by matrix-vector operations, we expect that the higher-order SB may be accelerated by parallel processing with, e.g., GPUs and FPGAs.

We thank Y. Sakai, T. Kashimata, K. Tatsumura, and R. Hidaka for helpful discussions.

\end{document}